\documentclass{birkjour}
\usepackage{amssymb}
\usepackage{amsmath}
\usepackage{ifthen}
\usepackage{graphicx}
\usepackage{amsfonts}
\usepackage{amscd}
\usepackage{amsxtra}
\usepackage{color}

 \newtheorem{thm}{Theorem}[section]
 \newtheorem{cor}[thm]{Corollary}
 \newtheorem{lem}[thm]{Lemma}
 
 \theoremstyle{definition}
 \newtheorem{defn}[thm]{Definition}
\newtheorem{rem}[thm]{Remark}
\newtheorem{conj}[thm]{Conjecture}

 \newtheorem{prob}[thm]{Problem}
 \numberwithin{equation}{section}

\begin{document}
%samy's command begins

%\usepackage{color}

\newcounter{alphabet}
\newcounter{tmp}
\newenvironment{Thm}[1][]{\refstepcounter{alphabet}%
\bigskip%
\noindent%
{\bf Theorem \Alph{alphabet}}%
\ifthenelse{\equal{#1}{}}{}{ (#1)}%
{\bf .} \itshape}{\vskip 8pt}
%\newcommand{\Ref}[1]{\setcounter{tmp}{\ref{#1}}\Alph{tmp}}

%% change begins by samy on 28-08-11
%\makeatletter
%%\newcommand{\Ref}[1]{\setcounter{tmp}{\ref{#1}}\Alph{tmp}}
%\newcommand{\Ref}[1]{\@ifundefined{r@#1}{}{\setcounter{tmp}{\ref{#1}}\Alph{tmp}}}
%\makeatother
%% change ends by samy on 28-08-11
%
%
%\newenvironment{Lem}[1][]{\refstepcounter{alphabet}%
%\bigskip%
%\noindent%
%{\bf Lemma \Alph{alphabet}}%
%{\bf .} \itshape}{\vskip 8pt}
%
%
%
%\newenvironment{Core}[1][]{\refstepcounter{alphabet}%
%\bigskip%
%\noindent%
%{\bf Corollary \Alph{alphabet}}%
%{\bf .} \itshape}{\vskip 8pt}

\newcommand{\co}{{\overline{\operatorname{co}}}}
\newcommand{\A}{{\mathcal A}}
\newcommand{\Bb}{{\mathbb B}}
\newcommand{\F}{{\mathcal F}}
\newcommand{\T}{{\mathbb T}}
\newcommand{\U}{{\mathcal U}}
\newcommand{\es}{{\mathcal S}}
\newcommand{\LU}{{\mathcal{LU}}}
\newcommand{\ZF}{{\mathcal{ZF}}}
\newcommand{\IR}{{\mathbb R}}
\newcommand{\IN}{{\mathbb N}}
\newcommand{\IC}{{\mathbb C}}
\newcommand{\ID}{{\mathbb D}}
\newcommand{\IB}{{\mathbb B}}
\newcommand{\K}{{\mathcal K}}
\newcommand{\X}{{\mathcal X}}
\newcommand{\PP}{{\mathcal P}}
\newcommand{\uhp}{{\mathbb H}}
\newcommand{\Z}{{\mathbb Z}}
\newcommand{\N}{{\mathcal N}}
\newcommand{\M}{{\mathcal M}}
\newcommand{\SCC}{{\mathcal{SCC}}}
\newcommand{\CC}{{\mathcal C}}
\newcommand{\st}{{\mathcal{SS}}}
\newcommand{\D}{{\mathbb D}}
\newcommand{\sphere}{{\widehat{\mathbb C}}}
\newcommand{\image}{{\operatorname{Im}\,}}
\newcommand{\Aut}{{\operatorname{Aut}}}
\newcommand{\real}{{\operatorname{Re}\,}}
\newcommand{\kernel}{{\operatorname{Ker}}}
\newcommand{\ord}{{\operatorname{ord}}}
\newcommand{\id}{{\operatorname{id}}}
\newcommand{\mob}{{\text{\rm M\"{o}b}}}
\newcommand{\Int}{{\operatorname{Int}\,}}
\newcommand{\Sign}{{\operatorname{Sign}}}
\newcommand{\diam}{{\operatorname{diam}}}
\newcommand{\inv}{^{-1}}
\newcommand{\area}{{\operatorname{Area}}}
\newcommand{\eit}{{e^{i\theta}}}
\newcommand{\esssup}{{\operatorname{esssup}}}
\newcommand{\dist}{{\operatorname{dist}}}
\newcommand{\arctanh}{{\operatorname{arctanh}}}
\newcommand{\ucv}{{\operatorname{UCV}}}

\def\be{\begin{equation}}
\def\ee{\end{equation}}
\newcommand{\sep}{\itemsep -0.01in}
\newcommand{\seps}{\itemsep -0.02in}
\newcommand{\sepss}{\itemsep -0.03in}
\newcommand{\bee}{\begin{enumerate}}
\newcommand{\eee}{\end{enumerate}}
\newcommand{\pays}{\!\!\!\!}
\newcommand{\pay}{\!\!\!}
\newcommand{\blem}{\begin{lem}}
\newcommand{\elem}{\end{lem}}
\newcommand{\bthm}{\begin{thm}}
\newcommand{\ethm}{\end{thm}}
\newcommand{\bcor}{\begin{cor}}
\newcommand{\ecor}{\end{cor}}
\newcommand{\beg}{\begin{example}}
\newcommand{\eeg}{\end{example}}
\newcommand{\begs}{\begin{examples}}
\newcommand{\eegs}{\end{examples}}
\newcommand{\bdefe}{\begin{defn}}
\newcommand{\edefe}{\end{defn}}
\newcommand{\bprob}{\begin{prob}}
\newcommand{\eprob}{\end{prob}}
\newcommand{\bques}{\begin{ques}}
\newcommand{\eques}{\end{ques}}
\newcommand{\bei}{\begin{itemize}}
\newcommand{\eei}{\end{itemize}}

\newcommand{\bde}{\begin{deter}}
\newcommand{\ede}{\end{deter}}
\newcommand{\bca}{\begin{case}}
\newcommand{\eca}{\end{case}}
\newcommand{\bcl}{\begin{claim}}
\newcommand{\ecl}{\end{claim}}
\newtheorem{claim}{Claim}
\renewcommand{\theclaim}{\!\!}
\newcommand{\bcon}{\begin{conj}}
\newcommand{\econ}{\end{conj}}
\newcommand{\bcons}{\begin{conjs}}
\newcommand{\econs}{\end{conjs}}
\newcommand{\bprop}{\begin{propo}}
\newcommand{\eprop}{\end{propo}}
\newcommand{\brs}{\begin{rems}}
\newcommand{\ers}{\end{rems}}
\newcommand{\bo}{\begin{obser}}
\newcommand{\eo}{\end{obser}}
\newcommand{\bos}{\begin{obsers}}
\newcommand{\eos}{\end{obsers}}
\newcommand{\bpf}{\begin{proof}}
\newcommand{\epf}{\end{proof}}
\newcommand{\ba}{\begin{array}}
\newcommand{\ea}{\end{array}}
\newcommand{\beq}{\begin{eqnarray}}
\newcommand{\beqq}{\begin{eqnarray*}}
\newcommand{\eeq}{\end{eqnarray}}
\newcommand{\eeqq}{\end{eqnarray*}}
\newcommand{\llra}{\longleftrightarrow}
\newcommand{\lra}{\longrightarrow}
\newcommand{\lla}{\longleftarrow}
\newcommand{\Llra}{\Longleftrightarrow}
\newcommand{\Lra}{\Longrightarrow}
\newcommand{\Lla}{\Longleftarrow}
\newcommand{\Ra}{\Rightarrow}
\newcommand{\La}{\Leftarrow}
\newcommand{\ra}{\rightarrow}
\newcommand{\la}{\leftarrow}
\newcommand{\ds}{\displaystyle}
\newcommand{\psubset}{\varsubsetneq}
%samy's command ends

%-------------------------------------------------------------------------
% editorial commands: to be inserted by the editorial office
%
%\firstpage{1}
%\volume{228}
%\DOI{003-0001}
%\seriesextraline{This is the Concrete Title of this Book\br H.E. R and S.T.C. W, Eds.}
%\firstpage{1}
%\title[Convolutions of harmonic mappings convex in one direction]
% {Landau-Bloch constants for functions in \\$\alpha$-Bloch spaces and Hardy spaces}

\title[Functions convex in some direction and Blaschke products]
{Geometric subfamily of functions convex in some direction and
Blaschke products}

%----------Author 1
\author[L. Li]{Liulan Li}
\address{College of Mathematics and Statistics\\
 (Hunan Provincial Key Laboratory of Intelligent Information Processing and Application),\\
Hengyang Normal University,\\
 Hengyang,  Hunan 421002,\\
People's Republic of China.}
\email{lanlimail2012@sina.cn}

%\thanks{The second author is currently on leave from the Indian Institute of Technology Madras, Chennai-600 036, India.}
%%----------Author 2

\author[S. Ponnusamy]{Saminthan Ponnusamy}
\address{Department of Mathematics,\\
Indian Institute of Technology Madras,\\
Chennai-600 036, India. }
\email{samy@iitm.ac.in}


\keywords{Univalent, starlike, convex and close-to-convex functions, harmonic functions, Blaschke products, John disk,
Schwarzian and pre-Schwarzian derivatives, quasiconformal mappings.\\ \\
{\tt The results of this article are a part of a report of June 2020, and an enlarged version of this
article was indeed with a journal in June 2020.}
%Correspondence should be addressed to S. Ponnusamy at {\tt samy@iitm.ac.in}
%\\ ${}^{\mathbf{*}}$ Corresponding author
}

\date{June 02, 2020}

\begin{abstract}
Consider the family of  locally univalent analytic functions $h$ in the unit disk $|z|<1$ with the normalization $h(0)=0$, $h'(0)=1$
and satisfying the condition
$${\real} \left( \frac{z h''(z)}{\alpha h'(z)}\right) <\frac{1}{2} ~\mbox{ for $z\in \ID$,}
$$
where $0<\alpha\leq1$. The aim of this article is to show that this family has several elegant properties such as involving Blaschke products,
Schwarzian derivative and univalent harmonic mappings.
\end{abstract}

%%% ----------------------------------------------------------------------
\maketitle

\section{Introduction}

For $r>0$, let $\mathbb{D}_r:=\{z\in\mathbb{C}:\,|z|<r\}$ and
$\mathbb{D}:=\mathbb{D}_1$, the open unit disk. Also,
$\mathbb{T}:=\partial \mathbb{D}=\{z\in\mathbb{C}:\,|z|=1\}$ denotes
the unit circle. Let ${\mathcal A}$ be the family of normalized
analytic functions $h$ in $\ID$ with the normalization $h(0)=0$, $h'(0)=1$
%\be\label{eq1}
%f(z)=z+\sum_{n=2}^{\infty}a_{n}z^{n}
%\ee
and ${\mathcal S}$ be the subfamily of ${\mathcal A}$ consisting of all univalent
functions in $\ID$. Denote by ${\mathcal S}^*$ the subfamily of ${\mathcal S}$
consisting of starlike functions,  i.e., $f(\ID)$ is a domain which is starlike with respect to the origin.
Note that ${\mathcal A}$ is a complete metric space.

In this article, we consider  the family $\mathcal{G}(\alpha)$ of locally univalent functions $h\in {\mathcal A}$ such that
\be\label{eq2}
{\real} \left( \frac{z h''(z)}{\alpha h'(z)}\right) <\frac{1}{2} ~\mbox{ for each $z\in \ID$,}
\ee
where $0<\alpha\leq1$. The formulation \eqref{eq2} clearly shows that these functions in $\mathcal{G}(\alpha)$ can be composed with a conformal
representation of the half-plane $\{w: \, {\real } w<1/2\}$ onto $\ID$. A candidate for such
function is $f^{-1}(w) = w/(1+w)$, where $f(z)= z/(1-z)$, $(f\in {\mathcal S})$. Recall that
$f$ plays the role of Koebe's function for convex domains.

Functions in the family $\mathcal{G}:=\mathcal{G}(1)$ are known to be univalent in $\D$ \cite{Oz} and thus,
$\mathcal{G}(\alpha) \subset \mathcal{G} \subset \mathcal{S}$ for $0<\alpha\leq1$.  In fact, it is
shown that $h'(z)\prec H_{\alpha}(z)=(1-z)^\alpha$  for $z\in \ID$ if $h\in \mathcal{G}(\alpha)$,
where $\prec$ is the usual subordination \cite{OPW-2013}.
Therefore, if $h\in \mathcal{G}$, then
$${\real} \left(h'(z)\right)>0 ~\mbox{ for $z\in \ID$},
$$
and thus $h$ is univalent in $\ID$.  Also, it is known
(see \cite[Equation (16)]{PoRa1}) that functions in ${\mathcal G}(\alpha )$ are starlike.
This family ${\mathcal G}(\alpha )$ has been investigated extensively. See  \cite{OPW-2013, OPW18, Um} and the references therein. It is trivial to prove that ${\mathcal G}(\alpha )$ is compact.
%is mentioned in Theorem \ref{compact}.

Moreover it is known that \cite{OPW-2013} if $f\in {\mathcal G}(\alpha )$ for some $0<\alpha \leq 1$ and $f(z)=\sum_{n=1}^{\infty}a_nz^n$, then
the following coefficient inequality holds
$$|a_n|\leq \frac{\alpha}{n(n-1)}   ~\mbox{ for $n\geq 2$},
$$
and equality is attained for the function $f_n$ such that $f_n'(z)=(1-z^{n-1})^{\alpha /(n-1)}$, $n\geq 2$.

%Recall that the logarithmic coefficients $\gamma_n$ of $f\in  {\es}$ are defined by the formula
%%\begin{equation}
%%\label{eq1-sec1}
%$$\frac{1}{2}\log\left (\frac{f(z)}{z} \right )=\sum_{n=1}^\infty \gamma_n(f)z^n \quad \mbox{ for } z\in \ID.
%$$
%The inequality $|\gamma_n| \le 1/n$ holds for $f\in\es^*$ but is not true for the full class $\es$, even in order
%of magnitude (see \cite[Theorem 8.4]{Du}). Related information and importance on the logarithmic coefficients $\gamma_n$ of $f\in {\es}$
%and some related geometric subfamilies may be obtained from \cite{dB85,MG86,PonShWir20,Roth07} and the references therein.
%Concerning the family $\mathcal{G}(\alpha)$, the following conjecture was proposed in \cite{PonShWir20}.
%
%\bcon
%{\it
%The logarithmic coefficients $\gamma_n$ of the functions in $\mathcal{G}(\alpha),\, \alpha\in (0,1],$ satisfy the inequalities
%$$|\gamma_n|\,\leq \,\frac{\alpha}{2n(n+1)},\quad n\in \IN.
%$$
%Equality is attained for $f'(z)\,=\,(1-z^n)^{\alpha/n}.$
%}
%\econ
%In \cite{PonShWir20}, among other things, the authors have verified this conjecture for $n=1,2,3$ but it remains open for $n\geq 4$.

%Moreover, functions in $\mathcal{G}(1)$ are proved to be starlike in
%$\D$; see Equation (16) in \cite{PR} and Theorem 1 in \cite{JO} (see
%also \cite{PV}). Thus, the family $\mathcal{G}(\alpha)$ is included in
%${\mathcal S}$ whenever $\alpha\in (0,1]$.

In this article the authors continue to investigate further properties about the family $\mathcal{G}(\alpha)$.
The article is organized as follows. In Section \ref{LP17-sec2}, we build material to characterize functions in
$\mathcal{G}(\alpha)$ in terms of finite Blaschke product (see Theorem \ref{with Blaschke produt}).
% and\ref{Blaschke produt}).
In Section \ref{LP17-sec3}, we obtain sharp estimate of Schwarzian derivative of functions
belonging to the family $\mathcal{G}(\alpha)$ (see Theorem \ref{Schwarzian derivative}). In Section \ref{LP17-sec4},
we use a characterization of functions from $\mathcal{G}(\alpha)$ and then establish univalent harmonic
mappings $f=h+\overline{g}$ with the analytic part $h$ from the family $\mathcal{G}(\alpha)$.

\section{Properties and characterization of functions in $\mathcal{G}(\alpha)$}\label{LP17-sec2}

Let $\Bb$ denote the set of the functions $\omega$ which are analytic
in $\D$ and satisfy $|\omega(z)|\leq 1$ for all $z\in\D$. Also,  consider the subfamily
$\Bb_0=\{\omega\in \Bb:\, \omega(0)=0 \}$.
% is convex and compact,
%which together with \cite[Theorem 7.12 ]{HM} gives the following.

\begin{thm}\label{changing real part} Suppose that
$h\in\mathcal{G}(\alpha)$ for some $ \alpha\in (0,1]$. Then
\be\label{eq3e} {\real} \left(\frac{z
h''(z)}{h'(z)}\right)
\leq\frac{\alpha}{2}-\frac{1}{2\alpha}(1-|z|^2)\left|\frac{h''(z)}{h'(z)}\right|^2 ~\mbox{ for all $z\in \ID$.}
\ee
Strict inequality holds for all $z\in\ID$ unless $h'(z)=(1-\zeta z)^\alpha$ for some $\zeta\in\T$.
\end{thm}

\bpf By assumption, \eqref{eq2} holds and thus, there exists an  $\omega \in \Bb_0$, i.e. $\omega:\,\ID\rightarrow \ID$ with $\omega(0)=0$, such that
$$\frac{z h''(z)}{h'(z)}=\alpha\frac{\omega(z)}{\omega(z)-1} ~\mbox{ for $z\in \ID$}.
$$
Since $\omega(0)=0$, we may set $\omega(z)=z\phi(z)$  so that the last equation reduces to
\be\label{eq4}
\frac{h''(z)}{h'(z)}=\alpha\frac{\phi(z)}{z\phi(z)-1}  ~\mbox{ for $z\in \ID$},
\ee
where $\phi\in \Bb$. Rewriting \eqref{eq4} yields
$$\phi(z)=\frac{h''(z)/h'(z)}{ (zh''(z)/h'(z))-\alpha} ~\mbox{ for $z\in \ID$}.
$$
Since $|\phi(z)|\leq 1$, the last equation yields that
$$\left|\frac{h''(z)}{h'(z)}\right|^2\leq \left|\frac{zh''(z)}{h'(z)}-\alpha\right|^2~\mbox{ for $z\in \ID$}.
$$
Expanding the last inequality cleary shows that \eqref{eq3e} holds.
%which together with the fact that $|\phi(z)|\leq 1$ implies
%$${\real} \left(\frac{z h''(z)}{h'(z)}\right)
%\leq\frac{\alpha}{2}-\frac{1}{2\alpha}(1-|z|^2)\left|\frac{h''(z)}{h'(z)}\right|^2 ~\mbox{ for $z\in \ID$}.
%$$
If the equality in \eqref{eq3e} holds for some point $z_0\in\ID$, then $|\phi(z_0)|=1$ and hence,
$\phi(z)\equiv \zeta$ for some $\zeta\in\T$. This gives by \eqref{eq4} that,
$$\frac{h''(z)}{h'(z)}=\frac{\alpha \zeta}{\zeta z -1}  ~\mbox{ for $z\in \ID$},
$$
which by integration shows that   $h'(z)=(1-\zeta z)^\alpha$ as desired.
 \epf

%\begin{thm}\label{boundary real part}
%Suppose that $h\in\mathcal{G}(\alpha)$ for some $ \alpha\in (0,1]$,
%and let $\phi\in \Bb$ that satisfies \eqref{eq4}. Then
%\be\label{eq5}
%{\real} \left(\frac{z h''(z)}{h'(z)}\right) = \frac{\alpha}{2} ~\mbox{ for all $z\in \T$}
%\ee
%if and only if  either $\phi(z)\equiv\zeta$ for some $\zeta\in\T$ or $\phi(z)$ is a finite Blaschke product.
%\end{thm}
%\bpf From \eqref{eq4}, it follows   that  \eqref{eq5} holds if and
%only if $|\phi(z)|=1$ for all $z\in \T$, from which the desired equivalent condition follows.
%%which is equivalent to that
%%$$\phi(z)\equiv\zeta\;\ \mbox{for}\;\ \mbox{some}\;\ \zeta\in\T$$
%%or $\phi(z)$ is a finite Blaschke product.
%\epf

\subsection{Characterization of functions from $\mathcal{G}(\alpha)$}
\begin{thm} \label{derivative}
Let $h\in\mathcal{G}(\alpha)$ for some $ \alpha\in (0,1]$. Then
\be\label{LP17-eq1}
h'(z)=\exp \left ( \alpha\int_{\T}\log(1-\zeta z)\,d\mu(\zeta ) \right )~\mbox{ for $z\in \ID$},
\ee
where $\mu$ is a probability measure on $\T$ so that $\int_{\T} d\mu(\zeta )=1$. Also, we have
$h'(z)\prec H_{\alpha}(z)$ for $z\in \ID$, where $H_{\alpha}(z)=(1-z)^\alpha$  and $\prec$ is the usual subordination \cite{Du,Pom}.
%\bee
%\item [{\rm (1)}]
%$$\frac{h''(z)}{h'(z)}=-\alpha\int_{\T} \frac{\zeta }{1-\zeta z}d\mu(\zeta ),$$
%
%\item [{\rm (2)}] $$\log h'(z)=\alpha\int_{\T}\log(1-\zeta z)d\mu(\zeta ),$$
%\eee
\end{thm}
\bpf Let $h\in\mathcal{G}(\alpha)$. Then, from the analytic characterization of the family given by \eqref{eq2}, we
have an equivalent condition
$${\real } \left( 1-\frac{2}{\alpha}\frac{z h''(z)}{h'(z)}\right) >0 ~\mbox{ for $z\in \ID$},
$$
and hence, by the Herglotz representation for analytic functions with positive real part in the unit disk, it
follows easily that
$$1-\frac{2}{\alpha}\frac{z h''(z)}{h'(z)}=\int_{\T}
\frac{1+\zeta z}{1-\zeta z}\,d\mu(\zeta ), ~\mbox{ i.e. }~\frac{h''(z)}{h'(z)}=-\alpha\int_{\T} \frac{\zeta }{1-\zeta z}\,d\mu(\zeta ), ~\mbox{ for $z\in \ID$},
$$
where $\mu$ is a probability measure on $\T$. Integrating the last relation gives
$$\log h'(z)=\alpha\int_{\T}\log(1-\zeta z)\,d\mu(\zeta )
$$
and the desired conclusion \eqref{LP17-eq1} follows if we perform exponentiation on both sides of this relation.

Next we show that $h'(z)\prec H_{\alpha}(z)=(1-z)^\alpha$. This is well-known and is used, for example, in  \cite{OPW-2013}.
Here is an alternate proof of this subordination property of the family $\mathcal{G}(\alpha)$. To do this, we
let $\varphi (z)=-\log(1-z)$. Then $\varphi(z)$ is univalent and convex in $\D$.
It is obvious that $\int_{\T}\log(1/(1-\zeta z))\,d\mu(\zeta )$ is a convex average
of values in $\varphi(\D)$. Therefore, there exists an analytic function
$\omega:\,\D\rightarrow\D$ with $\omega(0)=0$ such that
$$\int_{\T}\log(1-\zeta z)\,d\mu(\zeta )=\log(1-\omega(z)).
$$
This together with \eqref{LP17-eq1} imply that $h'(z)=\left(1-\omega(z)\right)^\alpha$ and hence,
$h'(z)\prec (1-z)^\alpha$  for $z\in \ID$.
\epf

\begin{thm}\label{with Blaschke produt}
Suppose that $h\in\mathcal{G}(\alpha)$ for some $ \alpha\in (0,1]$, and satisfies \eqref{eq4} for some $\phi\in \Bb$.
Then $\phi$ is a finite Blaschke product with degree $m\geq 1$ if and only if
%\be\label{eq6}
$$h'(z)=\prod^{m+1}_{k=1}(1-\zeta_k z)^{\alpha t_k},
$$
where $\zeta_k\in\T$ are distinct points, $0<t_k<1$ and $\sum^{m+1}_{k=1} t_k=1$.
\end{thm}
\bpf {\bf Necessity.}
First we observe that the family $\mathcal{G}(\alpha)$ is rotationally invariant in the sense that if
$h\in \mathcal{G}(\alpha)$ then for each $\theta \in [0,2\pi)$, $e^{-i\theta}h(e^{i\theta}z)$
belongs to $\mathcal{G}(\alpha)$. Therefore, if $\phi$ is a finite Blaschke product with degree
$m\geq 1$, then by rotating $h$, we may assume that
%\be\label{eq7}
$$\phi(z)=\prod^{m}_{k=1} \frac{z-b_k}{1-\overline{b_k} z},\;\
b_k\in\ID.
$$
Therefore we have
\be\label{eq8}
\frac{\phi(z)}{z\phi(z)-1}=\frac{\prod^{m}_{k=1}(z-b_k)}{z\prod^{m}_{k=1}(z-b_k)-\prod^{m}_{k=1}(1-\overline{b_k}
z)}.\ee
From the right side of \eqref{eq8}, we can see that $\frac{\phi(z)}{z\phi(z)-1}$ is a rational function with
poles at the roots of $z\phi(z)=1$. Since $z\phi(z)$ is a finite Blaschke product with degree
$m+1$, Theorem 3.4.10 in \cite{GMR} implies that there exist $m+1$ distinct roots $z_1,z_2, \ldots, z_{m+1}$ on $\T$
of $z\phi(z)=1$. These points are simple poles of $\frac{\phi(z)}{z\phi(z)-1}$. A partial fraction expansion
gives that
 \be\label{eq9}
\frac{\phi(z)}{z\phi(z)-1}=\sum^{m+1}_{k=1} \frac{t_k}{z-z_k},\ee
where $t_k\neq 0$ are complex constants. Since the right side of
\eqref{eq8} is the quotient of two monic polynomials with the
numerator of degree $m$ and the denominator of degree $m+1$,
\eqref{eq8} and \eqref{eq9} shows that $\sum^{m+1}_{k=1} t_k=1$.

By using \eqref{eq9}, we have
$$t_k=\lim_{z\rightarrow z_k} \frac{(z-z_k)\phi(z)}{z\phi(z)-1}=\lim_{z\rightarrow z_k} \frac{\phi(z)+(z-z_k)\phi'(z)}{\phi(z)+z\phi'(z)}
=\frac{\phi(z_k)}{\phi(z_k)+z_k\phi'(z_k)},$$
which yields that
$$t_k=\frac{1}{1+\frac{z_k\phi'(z_k)}{\phi(z_k)}}.
$$
By using  \cite[(3.4.7)]{GMR}, we have
$$\frac{z_k\phi'(z_k)}{\phi(z_k)}>0,
$$
which shows that $0<t_k<1$ . In view of the last observation, by using \eqref{eq4} and \eqref{eq9}, we easily have
\be\label{eq12}
\frac{h''(z)}{h'(z)}=\alpha\frac{\phi(z)}{z\phi(z)-1}
=\alpha\sum^{m+1}_{k=1}\frac{t_k }{z-z_k}=\alpha\sum^{m+1}_{k=1}\frac{t_k (-\overline{z_k})}{1-\overline{z_k}z},
\ee
which by integration yields that
$$\log h'(z)=\alpha\sum^{m+1}_{k=1} t_k\log(1-\zeta_k z), ~\mbox{ i.e., }~
h'(z)=\prod^{m+1}_{k=1}(1-\zeta_k z)^{\alpha t_k},
$$
where $\zeta_k=\overline{z_k}.$

{\bf Sufficiency.} Assume that $h'(z)=\prod^{m+1}_{k=1}(1-\zeta_k
z)^{\alpha t_k},$ where $\zeta_k\in\T$ are distinct points,
$0<t_k<1$ and $\sum^{m+1}_{k=1} t_k=1$. By \eqref{eq4}, a
calculation gives that
\be\label{eq13}
\alpha\frac{\phi(z)}{z\phi(z)-1}=\frac{h''(z)}{h'(z)}=\alpha\sum^{m+1}_{k=1}
\frac{t_k}{z-\overline{\zeta_k}} ~\mbox{ and }~ \phi(0)=\sum^{m+1}_{k=1}
t_k\zeta_k.\ee
Since $0<t_k<1$ and $\zeta_k\in\T$ are distinct
points, we have $\phi(0)\in\D$. By \cite[Theorem 3.5.2]{GMR}, it
suffices to prove that $\lim_{|z|\rightarrow 1^{-}} |\phi(z)|=1$. In order to prove this, we consider
the first equation \eqref{eq13} and obtain that
\be\label{eq14}
z\phi(z)=\frac{\sum^{m+1}_{k=1} \frac{t_k
z}{z-\overline{\zeta_k}}}{\sum^{m+1}_{k=1} \frac{t_k
z}{z-\overline{\zeta_k}}-1}.
\ee
%and in particular, we have $\phi(\overline{\zeta_k})=\zeta_k$ for $1\leq k\leq m+1.$
For $1\leq k\leq m+1$, we have
$$\lim_{z\rightarrow \overline{\zeta_k}} z\phi(z)=1.
$$
Moreover, for $|a|=1$, $a\neq \zeta$ and $|\zeta|=1$, one has
$${\real} \left( \frac{\zeta}{\zeta-\overline{a}}\right)={\real} \left(\frac{a\zeta}{a\zeta- 1}\right)=\frac{1}{2}
$$
and therefore,
$${\real} \left(\sum^{m+1}_{k=1} \frac{t_k \zeta}{\zeta-\overline{\zeta_k}}\right)=\frac{1}{2} \sum^{m+1}_{k=1} t_k =\frac{1}{2}\;\ \mbox{for}\;\
\mbox{all}\;\
\zeta\in\T\backslash \{\overline{\zeta_1},\overline{\zeta_1},\ldots,\overline{\zeta_{m+1}}\}.
$$
It follows from \eqref{eq14} that $\lim_{z\rightarrow \zeta}|z\phi(z)|=1$ for all
$\zeta\in\T \backslash \{\overline{\zeta_1},\overline{\zeta_1},\ldots,\overline{\zeta_{m+1}}\}$.
So $\lim_{|z|\rightarrow 1^{-}} |\phi(z)|=1$ and the proof is completed.
 \epf

%Theorems \ref{changing real part}, \ref{boundary real part} and \ref{with Blaschke produt} imply that
%
%\begin{thm}\label{Blaschke produt}
%Suppose that $h\in\mathcal{G}(\alpha)$ for some $ \alpha\in (0,1]$, and satisfies \eqref{eq4} for some $\phi\in \Bb$.
%Then $z\phi (z)$ is a finite Blaschke product with degree $m\geq 1$ if and only if
%%\be\label{eq15}
%$$h'(z)=\prod^{m}_{k=1}(1-\zeta_k z)^{\alpha t_k},
%$$
%where $\zeta_k\in\T$ are distinct points, $0<t_k\leq 1$ and
%$\sum^{m}_{k=1} t_k=1$.
%\end{thm}

\section{Pre-Schwarzian derivative and Schwarzian derivative}\label{LP17-sec3}

A series of results have been established by using relationship between the univalence
of a locally univalent analytic function and its Schwarzian derivative or pre-Schwarzian derivative.
The origin of such an approach is connected with the investigations
of Nehari \cite{N1} using Schwarzian derivative. Subsequently, this idea has been significantly developed by
a number of researchers.  See \cite[Section 8.5]{Du} and \cite{AvAk1975, CDO2011, Pom} for further detail.
%This approach can be elucidated in the simplest way in the following manner.

For a locally univalent analytic function $h$ in $\ID$, we define the
pre-Schwarzian derivative $P_h$ and the Schwarzian derivative $S_{h}$ by
$$P_{h}(z)=\frac{h''(z)}{h'(z)} $$
and
$$
 S_{h}(z):=P'_{h}(z)-\frac{1}{2}P^2_{h}(z)=\frac{d}{dz}\left (\frac{h''(z)}{h'(z)}\right )
-\frac{1}{2}\left (\frac{h''(z)}{h'(z)}\right )^2,
$$
respectively.
Note that $P_h$ can be derived from the Jacobian $J_h=|h'|^2$ of $h$, namely,
$$P_{h}(z)=\frac{\partial}{\partial z}(\log J_h).
$$
Their norms are defined by
$$\|T_h\|=\sup_{z\in\ID} (1-|z|^2)|P_{h}(z)| ~\mbox{ and }~ \|S_h\|=\sup_{z\in\ID}
(1-|z|^2)^2|S_{h}(z)|,
$$
respectively.
%For pre-Schwarzian derivative,  the Becker univalence criterion \cite{Bec72} is much deeper
%(cf. \cite{ Bec72, BePo84, DuShSh66}).

There are several well-known results which ensure that $f$ is univalent
in $\ID$ involving these two quantities, and, possibly, has an extension to $K$-quasiconformal mapping $\varphi$
of the extended complex plane $\IC_\infty=\IC \cup \{\infty\}$ onto itself, where $K=(1+k)/(1-k)$. That is, there exists
a quasiconformal homeomorphism $\varphi$ on $\IC_\infty$ such that $\varphi(z)=f(z)$ for $z\in \ID$ and $$\|\mu\|_\infty :=\esssup\{ |\mu (z)|:\, z\in \IC\} \leq k<1,$$
where $\mu =\varphi_{\overline{z}}/\varphi _z$.
For some historical and further discussion on these derivatives, we refer to \cite{AgraSahoo-2020,KimSu02,PonSahSuga14} and
the references therein.

Finally, we recall that a quasicircle in $\IC_\infty$ is the image of a circle under a quasiconformal mapping of the plane. A domain bounded by a quasicircle is called a quasidisk. We remark that every bounded quasidisk is known to be a John disk, but not the converse. For pre-Schwarzian derivative, the following Becker univalence criterion \cite{Bec72} is much deeper (cf. \cite{ Bec72, BePo84, DuShSh66}).

\begin{Thm}\label{quasiconformal extension}
{\rm \cite{Bec72, BePo84}}
If $\|T_h\|\leq  1$, then $h$ is univalent in $\ID$ and the constant $1$ is best possible. Moreover, if
$\|T_h\|\leq k<1$, then $h$  has a continuous extension $\widetilde{h}$ to $\overline{\ID}$ and $h(\T)$ is a quasicircle.
\end{Thm}

Indeed, Becker showed that,  if $\|T_h\|\leq k<1$ then $h$ has a $K$-quasiconformal extension to the whole complex plane $\IC$, where $K=(1+k)/(1-k)$.
In \cite{OPW-2013}, it is proved that $\|T_h\|\leq2\alpha$ for
$h\in \mathcal{G}(\alpha)$ and therefore, as a consequence of Becker's result, we have the following.

\begin{cor}\label{extension} Suppose that
$h\in\mathcal{G}(\alpha)$ for some $ \alpha\in (0,1]$. If $0<\alpha<\frac{1}{2}$, then $h(\T)$ is
a quasicircle and $h$ has a $(1+2\alpha)/(1-2\alpha)$-quasiconformal
extension to the whole complex plane $\IC$.
\end{cor}

For $\alpha\in [1/2,1)$, we have the following result.

\begin{thm}\label{quasidisk}
Suppose that $h\in\mathcal{G}(\alpha)$ for some $ \alpha\in (0,1]$. Then $h(\ID)$ is a quasidisk.
\end{thm}
\bpf Let $h\in\mathcal{G}(\alpha)$. Then  $h'(z)\prec H_{\alpha}(z)$, where $H_{\alpha}(z)=(1-z)^\alpha$.
Then $D=H_{\alpha}(\ID)$ is a bounded domain which is contained in the right half-plane
$${\mathbb H}=\{w:\ {\real }w>0\}.
$$
Let $K$ be a closed disk with center in ${\mathbb H}\setminus\overline{D}$.
Then  $K$ is a compact subset and for each $r>0$, one case
$$r K\bigcap ({\mathbb H}\setminus D)\neq \emptyset .
$$
Moreover, as $h'(z)\prec H_{\alpha}(z)$, we have $h'(\ID)\subset D$.
This fact and the result in \cite{CG} imply that $h(\ID)$ is a quasidisk.
\epf

%Let $D$ be the image of $\ID$ under $(1-z)^\alpha$. Then it is
%a simple exercise to see that $D$ is bounded by the sinusoidal
%spiral $\Omega (\alpha)$, where \beqq \Omega (\alpha)&=& \left\{\rho
%e^{i\theta}:\, \rho=\left(2\cos\frac{\theta}{\alpha}\right)^\alpha,~
%\mbox{ and }
%-\frac{\alpha\pi}{2}<\theta\leq \frac{\alpha\pi}{2}\right\}\\
%&=&\left\{w:\, {\rm Re\,}w >0 ~\mbox{ and } {\rm Re\,}\big (w^{-1/\alpha}\big) \, =\frac{1}{2} \right \},
%\eeqq
%which is clearly contained in the right half-plane $R=\{w:\ {\real }w>0\}$. Let $K$ be the closed
%disk with center $4$ and radius $1$. Then $K\subset R$ is a compact
%subset and $r K\bigcap (R\setminus D)\neq \emptyset$ for all $r>0$.
%Now we let $h\in\mathcal{G}(\alpha)$. Then, by Lemma \ref{derivative}, $h'(z)\prec (1-z)^\alpha$ for $z\in \ID$.
%This fact and the result in \cite{CG} imply that $h(\ID)$ is a quasidisk.
%\epf

Note also that for $\alpha \in (0,1)$,  $h\in \mathcal{G}(\alpha)$ implies that $\|T_h\|\leq 2\alpha<2$ and hence, $f(\ID)$ is a
John disk.

Next we recall the following well-known results from \cite{N1} (see also \cite{N2,N3} which deals with
general situation). As shown by Hille \cite{Hille}, the constant $2$ below is best possible.

\begin{Thm}\label{Nehari}
If $f\in {\mathcal S}$, then we have the sharp inequality $\|S_h\|\leq  6$ and the number $6$ is best possible.
Conversely, if $\|S_h\|\leq  2$, then $h$ is univalent in $\ID$ and the number $2$ is best possible.
Moreover, if $\|S_h\|\leq k<2$, then $h$  has a continuous extension to the whole complex plane.
\end{Thm}

\begin{lem}\label{LP17-lem-extra1}
Let $h\in\mathcal{G}(\alpha)$ for some $\alpha \in (0,1)$ such that $h'(z)=(1-\zeta z)^\alpha$ for all $z\in \ID$ and
$\zeta\in\T$. Then $\|S_h\|=2\alpha(2+\alpha)$.
\end{lem}
\bpf Proof follows from a direct computation. Indeed, for the function $h'(z)=(1-\zeta z)^\alpha$, we have
$$S_h(z) =-\frac{\alpha \zeta ^2}{(1-\zeta z)^2} -\frac{1}{2} \frac{\alpha ^2\zeta ^2}{(1-\zeta z)^2} =-\frac{\alpha(2+\alpha)}{2}
\frac{\zeta ^2}{(1-\zeta z)^2}
$$
showing that $\|S_h\|=2\alpha(2+\alpha).$
\epf

%Naturally, one may ask
%\bprob
%What is the bound for the Schwarzian derivative of $h\in\mathcal{G}(\alpha)$?
%\eprob
%
%Our next result answers this question.

In the following, we determine the bound for the Schwarzian derivative of $h\in\mathcal{G}(\alpha)$.

\begin{thm}\label{Schwarzian derivative}
Suppose that $h\in\mathcal{G}(\alpha)$ for some $ \alpha\in (0,1]$. Then $\|S_h\|\leq2\alpha(2+\alpha),$ where
the equality is attained by $h(z)$ which is obtained from $h'(z)=(1-\zeta z)^\alpha$ for some $\zeta\in\T$.
\end{thm}
\bpf Let $h\in\mathcal{G}(\alpha)$. Then, by using \eqref{eq4}, we have
$$\frac{h''(z)}{h'(z)}=\alpha\frac{\phi(z)}{z\phi(z)-1},
$$
for some $\phi\in \Bb$. Therefore, by Schwarz-Pick inequality, we have
\be\label{S-Pineq-eq1}
(1-|z|^2)|\phi'(z)| \leq  1-|\phi(z)|^2 ~\mbox{ for $z\in \ID$}.
\ee
Now, a calculation yields
$$S_h(z)= -\alpha\left (\frac{2\phi'(z)+(2+\alpha)\phi^2(z)}{2(z\phi(z)-1)^2}\right )
$$
so that
\beqq
%S_h(z)&=& -\alpha\frac{2\phi'(z)+(2+\alpha)\phi^2(z)}{2(z\phi(z)-1)^2}\\
\frac{1}{\alpha}|S_h|(1-|z|^2)^2&\leq&\frac{|\phi'(z)|(1-|z|^2)^2}{|z\phi(z)-1|^2}+\frac{(2+\alpha)|\phi(z)|^2(1-|z|^2)^2}{2|z\phi(z)-1|^2},
\eeqq
which, by using  the Schwarz-Pick inequality \eqref{S-Pineq-eq1}, reduces to
%\be\label{eq16}
\be\label{S-Pineq-eq2}
\frac{1}{\alpha}|S_h|(1-|z|^2)^2\leq\frac{(1-|\phi(z)|^2)(1-|z|^2)}{|z\phi(z)-1|^2}+\frac{(2+\alpha)|\phi(z)|^2(1-|z|^2)^2}{2|z\phi(z)-1|^2}.
\ee

If there exists a $z_0\in\ID$ such that $|\phi(z_0)|=1$, then $\phi(z)\equiv\zeta$ for some $\zeta\in\T$ and $h'(z)=(1-\zeta z)^\alpha$.  By Lemma \ref{LP17-lem-extra1}, we obtain that
%Therefore,
%$$S_h(z) =-\frac{\alpha \zeta ^2}{(1-\zeta z)^2}
%-\frac{1}{2} \frac{\alpha ^2\zeta ^2}{(1-\zeta z)^2} =-\frac{\alpha(2+\alpha)}{2} \frac{\zeta ^2}{(1-\zeta z)^2}
%$$
%showing that
$$\|S_h\|=2\alpha(2+\alpha).
$$

If $|\phi(z)|<1$ for all $z\in\ID$, then, by the inequality \eqref{S-Pineq-eq2}, we only need to prove that
$$\frac{(1-|\phi(z)|^2)(1-|z|^2)}{|z\phi(z)-1|^2}+\frac{(2+\alpha)|\phi(z)|^2(1-|z|^2)^2}{2|z\phi(z)-1|^2}\leq 2(2+\alpha),
$$
%$$2\frac{(1-|\phi(z)|^2)(1-|z|^2)}{|z\phi(z)-1|^2}+\frac{(2+\alpha)|\phi(z)|^2(1-|z|^2)^2}{|z\phi(z)-1|^2}\leq 4(2+\alpha),$$
or equivalently,
$$2(1-|z|^2)-2(1-|z|^2)|\phi(z)|^2+(2+\alpha)(1-2|z|^2+|z|^4)|\phi(z)|^2
$$
$$\leq 4(2+\alpha)\left(1-2{\real}(z\phi(z))+|z|^2|\phi(z)|^2\right).
$$
Simplifying the last inequality, one obtains an equivalent inequality
$M(\alpha) \leq0,$ where
\begin{align*}
%\nonumber
M(\alpha)& =\left((2+\alpha)|z|^4-(10+6\alpha)|z|^2+\alpha\right)|\phi(z)|^2+8(2+\alpha){\real}(z\phi(z))\\
& \hspace{2cm} -2\left(|z|^2+2\alpha+3\right).
\end{align*}
\begin{claim}%\label{c1}
$M(\alpha)$ is decreasing with respect to $\alpha$ for $\alpha\in [0, 1]$.
\end{claim}

By assumption, we have
\beqq
M'(\alpha)&=&\left(|z|^4-6|z|^2+1\right)|\phi(z)|^2+8{\real}(z\phi(z))-4\\
&=&\left(|z|^4-2|z|^2+1\right)|\phi(z)|^2-4|\phi(z)|^2+8{\real}(z\phi(z))-4\\
&=&(1-|z|^2)^2|\phi(z)|^2-4\left|z\phi(z)-1\right|^2\\
&=&\left((1-|z|^2)|\phi(z)|-2|z\phi(z)-1| \right)\left((1-|z|^2)|\phi(z)|+2|z\phi(z)-1| \right).
\eeqq
Now, we let $N(\alpha)=(1-|z|^2)|\phi(z)|-2|z\phi(z)-1|$. Then
\beqq
N(\alpha)&\leq& (1-|z|^2)|\phi(z)|-2(1-|z\phi(z)|)\\
&=&  (1-|z|^2+2|z|)|\phi(z)|-2\\
&=&  -(1-|z|)^2|\phi(z)|-2(1-|\phi(z)|)<0,\;\ z\in\ID. \eeqq
The above estimate together with the expression of $M'(\alpha)$ yields
$$M'(\alpha)<0,\;\ z\in\ID
$$
and  the claim is proved.

Consequently,  $M(\alpha)\leq M(0)$ for $\alpha\in (0,1]$. To complete the proof
it suffices to show that $M(0)\leq 0$ for all $z\in\ID$. Now,
%\begin{claim}\label{c2}
%$M(0)\leq 0$ for all $z\in\ID$.
%\end{claim}
$$M(0)=\left(2|z|^4-10|z|^2\right)|\phi(z)|^2+16{\real}(z\phi(z))-2\left(|z|^2+3\right)
$$
so that
\beqq
M(0)/2&=&\left(|z|^4-|z|^2\right)|\phi(z)|^2-4|z|^2|\phi(z)|^2+8{\real}(z\phi(z))-|z|^2-3\\
&=& -(1-|z|^2)|z\phi(z)|^2-4|z\phi(z)-1|^2+(1-|z|^2)\\
&=& (1-|z|^2)(1-|z\phi(z)|^2)-4|z\phi(z)-1|^2\\
&<& 4(1-|z|)(1-|z\phi(z)|)-4|z\phi(z)-1|^2\\
%&\leq&4[(1-|z|)|1-z\phi(z)|-|1-z\phi(z)|^2]\\
&<&4|1-z\phi(z)|\left((1-|z|)-|1-z\phi(z)|\right)<0
\eeqq
%By using the two claims, we have
and the proof of the theorem is complete.
\epf

\section{Univalent harmonic mappings with analytic part from $\mathcal{G}(\alpha)$}\label{LP17-sec4}

We consider the family of complex-valued functions $f=u+iv$ defined in the unit disk $\ID$, where $u$ and $v$ are real harmonic in $\ID$. Such
functions can be represented as $f=h+\overline{g}$, where $h$ and $g$ are analytic in $\ID$.
This representation is unique up to an additive constant and thus, without loss of generality we may assume
that $g(0)=0$.  The Jacobian $J_{f}$ of $f$ is given by $J_{f}(z)=|h'(z)|^{2}-|g'(z)|^{2}$, and thus,
a locally univalent function $f$ is sense-preserving if $J_{f}(z)>0$ in $\ID$. Consequently,
a harmonic mapping $f$ is locally univalent and sense-preserving in $\ID$ if and only if $J_{f}(z)>0$ in $\ID$ (cf. \cite{L1936}); or equivalently if $h'\neq 0$ in $\ID$ and the dilatation $\omega_{f}$ satisfies the
Beltrami equation of the second kind, $\overline{f_{\overline{z}}}=\omega _f f_{z}$, where $\omega_f$ is an analytic
function with $|\omega_f(z)|<1$ for $z\in \ID$ \cite{Clunie-Small-84, Du2}.
Note that $\omega_f =g'/h'$.

If a locally univalent and sense-preserving harmonic mapping $f=h+\overline{g}$ on $\ID$ satisfies the condition
$ |\omega_{f}(z)| \leq k<1$ for $\ID$,
then $f$ is called $K$-quasiregular harmonic mapping in $\ID$, where $K=\frac{1+k}{1-k}\geq 1$. We refer to \cite{Clunie-Small-84}
for several properties of univalent harmonic mappings together with its various subfamilies. In particular, here is a sufficient condition
for close-to-convexity of harmonic mappings due to Clunie and Sheil-Small \cite{Clunie-Small-84}. See also \cite{Du2}.

Recall that a domain $D$ is linearly connected if there exists a positive constant $M<\infty$ such
that any two points $z, w\in D$ are joined by a path $\gamma\subset D$ of length (cf. \cite{CG})
$$\ell(\gamma)\leq M|z-w| ,
$$
or equivalently $\diam (\gamma)\leq M|z-w|$. We point out that a bounded linearly connected domain is
a Jordan domain, and for piecewise smoothly bounded domains, linear connectivity is equivalent to the boundary's having no inward-pointing cusps. Also, if this inequality holds with $M=1$, then the linearly connected domain
is convex.

\begin{thm}\label{univalent harmonic functions}
Suppose that $h\in\mathcal{G}(\alpha)$ for some $ \alpha\in (0,1]$. Then there exists $c > 0$ such that
every harmonic mapping $f=h+\overline{g}$ with dilatation $|\omega(z)| < c$,
is univalent in $\ID$, where the constant $c$ depends only on the domain $h(\ID)$.
\end{thm}
\bpf Since $h(\ID)$ is a quasidisk by Theorem \ref{quasidisk}, $D=h(\ID)$ is linearly connected. The conclusion follows
from \cite[Theorem 1]{CR}.
\epf

\begin{thm}\label{univalent harmonic functions 2}
Let $h\in\mathcal{G}(\alpha)$ for some $ \alpha\in (0,1/2)$, and $f=h+\overline{g}$ be a sense-preserving harmonic
mapping with dilatation $\omega=g'/h'$. If
\be\label{eq17}
|\omega(z)|\leq1-\alpha|z|(1+|z|) ~\mbox{ for all $z\in\ID$,}
\ee
then $f$ is univalent in $\ID$.
\end{thm}
\bpf Since $h\in\mathcal{G}(\alpha)$, we have $h'(z)\prec
(1-z)^\alpha$ by Theorem \ref{derivative} and hence $h'(z)\neq0$ for
$z\in \ID$.  Recall Equation \eqref{eq4}:
$$\frac{h''(z)}{h'(z)}=\alpha\frac{\phi(z)}{z\phi(z)-1}~\mbox{ for all $z\in\ID$,}
$$
where $\phi\in \Bb$. By computation, we obtain that
$$ (1-|z|^2)\left |\frac{zh''(z)}{h'(z)}\right |\leq \alpha|z|(1+|z|)~\mbox{ for all $z\in\ID$,}
$$
which together with the result in \cite{FNS} implies that $f$ is univalent in $\ID$.
 \epf

\begin{rem}%\br
Note that \eqref{eq17} obviously holds whenever the dilatation satisfies the condition $|\omega(z)|\leq 1-2\alpha$ in $\ID$, where $\alpha \in (0,1/2)$.
\end{rem}%\er

\subsection*{Funding.}
The work of Liulan Li is supported by Hunan Provincial Natural Science Foundation of China (No. 2024JJ7063),
the Science and Technology Plan Project of
Hunan Province (No. 2016TP1020), and the
Application-Oriented Characterized Disciplines, Double First-Class
University Project of Hunan Province (Xiangjiaotong [2018]469).

%The work of the first author is supported by NSF of China (No.
%11571216), the Applied Characteristic Discipline Program in Hunan
%Province, the Science and Technology Plan Project of Hunan Province
%(No. 2016TP1020) and the Science and Technology Plan Project of
%Hengyang City (2017KJ183).

\subsection*{Availability of data and material.} Not applicable

\subsection*{Code availability.} Not applicable

\subsection*{Conflicts of interest/Competing interests.} There is no competing interests.

\end{document}